\documentclass[11pt]{article}
\usepackage{amsfonts}
\usepackage{color}
\newtheorem{theo}{Theorem}[section]
\newtheorem{lem}[theo]{Lemma}

\newcommand{\mysection}[1]{\section{#1} \setcounter{equation}{0}}
\newcommand{\proof}{{\sc Proof.} \quad}
\newcommand{\proofc}{{\sc Proof} \ }
\newcommand{\be}{\begin{equation} \label}
\newcommand{\ee}{\end{equation}}
\newcommand{\bea}{\begin{eqnarray}\label}
\newcommand{\eea}{\end{eqnarray}}
\newcommand{\bas}{\begin{eqnarray*}}
\newcommand{\eas}{\end{eqnarray*}}
\newcommand{\bit}{\begin{itemize}}
\newcommand{\eit}{\end{itemize}}
\newcommand{\qed}{\hfill$\Box$ \vskip.2cm}
\newcommand{\nn}{\nonumber}
\newcommand{\R}{\mathbb{R}}

\newcommand{\pO}{\partial\Omega}

\newcommand{\hra}{\hookrightarrow}
\newcommand{\io}{\int_\Omega}
\newcommand{\bom}{\overline{\Omega}}
\newcommand{\abs}{\\[5pt]}

\newcommand{\tm}{T_{max}}

\newcommand{\proj}{{\cal P}}
\newcommand{\ks}{K_S}

\newcommand{\oy}{\overline{y}}
\newcommand{\ts}{t_\star}
\newcommand{\hc}{\widehat{c}}
\newcommand{\hu}{\widehat{u}}
\newcommand{\I}{I_p}
\begin{document}
\enlargethispage{10mm}
\title{Does fluid interaction affect regularity in the three-dimensional Keller-Segel system with saturated
sensitivity?}
\author{
Michael Winkler\footnote{michael.winkler@math.uni-paderborn.de}\\
{\small Institut f\"ur Mathematik, Universit\"at Paderborn,}\\
{\small 33098 Paderborn, Germany} }
\date{}
\maketitle
\begin{abstract}
\noindent 
  A class of Keller-Segel-Stokes systems generalizing the prototype
  \bas
    	\left\{ \begin{array}{rcl}
    	n_t + u\cdot\nabla n &=& \Delta n  - \nabla \cdot \Big(n(n+1)^{-\alpha}\nabla c\Big), \\[1mm]
    	c_t + u\cdot\nabla c &=& \Delta c-c+n, \\[1mm]
     	u_t +\nabla P  &=& \Delta u + n \nabla \phi + f(x,t), \qquad \nabla\cdot u =0, 
    	\end{array} \right.
	\qquad \qquad (\star)
  \eas
  is considered in a bounded domain $\Omega\subset\R^3$, where $\phi$ and $f$ are given sufficiently
  smooth functions such that $f$ is bounded in $\Omega\times (0,\infty)$.\abs
  It is shown that under the condition that		
  \bas
	\alpha>\frac{1}{3},
  \eas
  for all sufficiently regular initial data
  a corresponding Neumann-Neumann-Dirichlet initial-boundary value
  problem possesses a global bounded classical solution.
  This extends previous findings asserting a similar conclusion only under the stronger assumption 
  $\alpha>\frac{1}{2}$.\abs
  In view of known results on the existence of exploding solutions when $\alpha<\frac{1}{3}$,
  this indicates that with regard to the occurrence of blow-up the criticality of the decay rate $\frac{1}{3}$,
  as previously found for the fluid-free counterpart of ($\star$), 
  remains essentially unaffected by fluid interaction of the type considered here.\abs
\end{abstract}
\newpage
\section{Introduction}\label{intro}
One of the most characteristic mathematical features of the classical Keller-Segel system, in its simplest form given by
\be{KS}
	\left\{ \begin{array}{l}
	n_t=\Delta n - \nabla \cdot (n\nabla c), \\[1mm]
	c_t=\Delta c-c+n,
	\end{array} \right.
\ee
consists in its ability to generate singular behavior by enforcing finite-time blow-up of some solutions 
in spatially two- or higher-dimensional situations (\cite{herrero_velazquez}, \cite{win_JMPA}). 
Well-established as a model for the collective behavior in populations of cells chemotactically biased 
by a signal substance produced by themselves,
(\ref{KS}) thus may well describe phenomena of spontaneous cell aggregation arising in various experimental contexts
(\cite{hillen_painter2009}). 
In order to adequately describe chemotactic migration also in biological frameworks in which such an emergence 
of unbounded population densities seems unrealistic, considerable efforts have been undertaken since the 
introduction of (\ref{KS}) (\cite{KS70}) to develop modified variants thereof
in which the occurrence of explosions is a priori ruled out. \abs
One frequently discussed and in its mathematical consequences quite comprehensively understood direction of refinement
consists in assuming the cell motility to depend differently on the population density than supposed in (\ref{KS}),
especially at large densities; 
this may lead to certain saturation effects in the cross-diffusion term, or to 
nonlinear diffusivities e.g.~in the sense of a porous medium-type enhancement of diffusion at large densities,
or to a combination of both 
(see e.g.~the survey \cite{hillen_painter2009}).
Focusing here on the former type of modification, as reflected in the variant
\be{KSM}
	\left\{ \begin{array}{l}
	n_t=\Delta n - \nabla \cdot (nS(n)\nabla c), \\[1mm]
	c_t=\Delta c-c+n,
	\end{array} \right.
\ee
of (\ref{KS}) with nonnegative $S(n)$ possibly becoming small at large values of $n$, we may interpret the
corresponding literature as identifying the decay rate of the prototypical choice
\be{proto}
	S(n)=(n+1)^{-\frac{N-2}{N}},
	\qquad n\ge 0,
\ee
as critical for the occurrence of blow-up in the spatially $N$-dimensional version of (\ref{KSM}): 
Indeed, if $N\ge 2$ and $S\in C^2([0,\infty))$ is such that 
\be{S_sub}
	S(n) \le \ks (n+1)^{-\alpha} 
	\qquad \mbox{for all } n\ge 0
\ee
and some
$\ks>0$ and $\alpha>\frac{N-2}{N}$, then for all reasonably regular nonnegative initial data 
the no-flux initial-boundary value problem for (\ref{KSM}) in smoothly bounded domains $\Omega\subset \R^N$
possesses a globally defined bounded classical solution (\cite{horstmann_win}, \cite{taowin_subcrit});
on the other hand, if 
\be{S_super}
	S(n) \ge \ks'(n+1)^{-\alpha'}
	\qquad \mbox{for all } n\ge 0
\ee
and some $\ks'>0$ and $\alpha'<\frac{N-2}{N}$,
then in each ball $\Omega\subset \R^N$ there exist solutions which become unbounded (\cite{cieslak_stinner_JDE2015},
\cite{win_collapse}).\abs
It is the purpose of the present work to study the question how far this borderline role of the behavior (\ref{proto})
may be affected by interaction of cells with a liquid environment, where intending to incorporate an assumption 
underlying the model development in \cite{goldstein2005} we will suppose that this interaction occurs not only through
transport but possibly also through a buoyancy-driven feedback of cells to the fluid velocity.
Indeed, numerical evidence suggests that the combination of these mechanisms may at least enforce a 
delay in blow-up of some solutions to an accordingly modified two-dimensional variant of (\ref{KS})
(\cite{lorz2012}).
More drastically, a recent rigorous analytical result shows that 
even in absence of any influence of cells on the fluid motion, a purely transport-determined interplay in fact may
fully suppress blow-up in the sense that for widely arbitrary fixed 
initial data one can construct a solenoidal fluid velocity
field such that a corresponding initial value problem associated with an either two- or three-dimensional variant
of (\ref{KS}) possesses globally bounded solutions (\cite{kiselev_ARMA2016}).\abs
With our focus slightly differing from that in the latter study, we will henceforth concentrate on the problem of
deciding	
whether for some given sensitivity parameter function $S$, in the extension 
\be{0}
    	\left\{ \begin{array}{rcll}
    	n_t + u\cdot\nabla n &=& \Delta n  - \nabla \cdot \Big(nS(n)\nabla c\Big), 
	\qquad & x\in\Omega, \ t>0,\\[1mm]
    	c_t + u\cdot\nabla c &=& \Delta c-c+n, 
	\qquad & x\in\Omega, \ t>0,   \\[1mm]
     	u_t +\nabla P  &=& \Delta u + n \nabla \phi + f(x,t), \qquad \nabla\cdot u =0, 
	\qquad & x\in\Omega, \ t>0,  \\[1mm]
	& & \hspace*{-28mm}
	\frac{\partial n}{\partial\nu}=0, \quad
	\frac{\partial c}{\partial\nu}=0, \quad
	u=0,
	\qquad & x\in\pO, \ t>0, \\[1mm]
	& & \hspace*{-28mm}
	n(x,0)=n_0(x), \quad c(x,0)=c_0(x), \quad u(x,0)=u_0(x), 
	\qquad & x\in\Omega,
    	\end{array} \right.
\ee
of the no-flux initial-boundary value problem for (\ref{KSM})
it is {\em at all} possible to observe the occurrence of blow-up for {\em some} solution in presence
of {\em some} suitably regular gravitational potential $\phi$ and external fluid force $f$ in {\em some} bounded
domain $\Omega\subset\R^N$.\abs
Within this problem setting it then immediately becomes clear on letting $u_0\equiv 0$, $\phi\equiv 0$, $f\equiv 0$
and $\Omega:=B_1(0)\subset\R^N$
that assuming (\ref{S_super}) to be valid for some $\ks'>0$ and $\alpha'<\frac{N-2}{N}$ trivially remains
sufficient for the existence of some exploding solutions in (\ref{0}) as well. 
In the case $N=2$, this condition in fact appears to stay essentially optimal also for (\ref{0})
in view of recent results asserting global existence of bounded classical solutions for all suitably regular
initial data at least when $f\equiv 0$, thus ruling out any blow-up phenomenon (\cite{wang_xiang_JDE2015}),
even in the more complicated case when the fluid flow is governed by an associated version of the full
Navier-Stokes equations (\cite{wang_win_xiang1}).\abs
In the three-dimensional version of (\ref{0}), the seemingly only available result on global existence and boundedness of
classical solutions for arbitrarily large initial data relies on the requirement that (\ref{S_sub}) holds for
some $\ks>0$ and $\alpha>\frac{1}{2}$ (\cite{wang_xiang_JDE2016}), thus leaving open the question how far the 
value $\frac{1}{3}$ accordingly appearing in (\ref{proto}) 
continues to play the role of a critical blow-up exponent for (\ref{0});
after all, under the mere assumption that (\ref{S_sub}) be valid with some $\ks>0$ and $\alpha>\frac{1}{3}$, 
certain global generalized solutions could be constructed for the actually even more complex  
Keller-Segel-Navier-Stokes variant of (\ref{0}) (\cite{wangyulan_3d_nasto}, cf.~also \cite{yifu_wang_JDE2017}), 
but unless in cases when
suitable additional smallness conditions on the initial data are imposed (\cite{kozono}) the knowledge
on their boundedness features is yet quite poor.\abs
{\bf Main results: Criticality of the decay exponent $\frac{1}{3}$.} \quad
The main outcome of this study reveals that the validity of (\ref{S}) with some $\ks>0$ and $\alpha>\frac{1}{3}$
is actually sufficient to exclude any singularity formation also in the full chemotaxis-Stokes system (\ref{0}) 
unde reasonable assumptions on $\phi, f$ and the initial data,
thereby indicating, in the sense specified above, that the possibility of observing blow-up in a suitable 
constellation remains unaffected by fluid interaction of the considered type.\abs
To make this more precise, let us consider (\ref{0})
in a bounded domain $\Omega \subset \R^3$, where for simplicity we shall assume that
\be{phi_f}
	\phi\in C^2(\bom)
	\qquad \mbox{and} \qquad
	f\in C^1(\bom\times [0,\infty);\R^3) \cap L^\infty(\Omega\times (0,\infty);\R^3),
\ee
and where we shall suppose throughout the sequel that
$S\in C^2([0,\infty))$ satisfies
\be{S}
    	|S(n)| \le \ks (n+1)^{-\alpha}
	\qquad \mbox{for all } n \ge 0
\ee
with some $\alpha>0$ and $\ks>0$.
The initial data in (\ref{0}) will be assumed to be such that
\be{init}
    	\left\{	\begin{array}{l}
    	n_0 \in C^0(\bar\Omega) \quad \mbox{ with } n_0\ge 0, \\
    	c_0 \in W^{1,\infty}(\Omega) \quad \mbox{ with $c_0\ge 0$ 
		\quad and}\\
    	u_0 \in D(A^{\beta}) \quad \mbox{for some $\beta\in (\frac{3}{4},1)$,}
    	\end{array} \right.
\ee
where $A=-\proj \Delta$ represents the Stokes operator in 
$L^2_\sigma(\Omega):=\{ \varphi\in L^2(\Omega;\R^3) \ | \ \nabla \cdot \varphi=0 \}$,
with its domain given by $D(A):=W^{2,2}(\Omega;\R^3)\cap W_0^{1,2}(\Omega;\R^3) \cap L^2_\sigma(\Omega)$,
and with $\proj$ denoting the Helmholtz projection from $L^2(\Omega;\R^3)$ into $L^2_\sigma(\Omega)$.\abs
In this context, our main results read as follows.
\begin{theo}\label{theo14}
  Let $\Omega\subset\R^3$ be a bounded domain with smooth boundary, let $\phi$ and $f$ satisfy (\ref{phi_f}),
  and let $S\in C^2([0,\infty))$ be such that (\ref{S}) holds with some
  \be{alpha}
	\alpha>\frac{1}{3}.
  \ee
  Then for all $n_0, c_0$ and $u_0$ fulfilling (\ref{init}), the problem (\ref{0}) possesses a global
  classical solution $(n,c,u,P)$, uniquely determined by the inclusions
  \be{14.1}
	\left\{ \begin{array}{l}
	n \in C^0(\bom\times [0,\infty)) \cap C^{2,1}(\bom\times (0,\infty)), \\
	c \in \bigcap_{p>3} C^0([0,\infty);W^{1,p}(\Omega)) \cap C^{2,1}(\bom\times (0,\infty)), \\
	u \in C^0([0,\infty);D(A^\beta)) \cap C^{2,1}(\bom\times [0,\infty);\R^3), \\
	P \in C^{1,0}(\bom\times (0,\infty)),
	\end{array} \right.
  \ee
  for which
  $n \ge 0$ and $c \ge 0$ in $\Omega\times (0,\infty)$.
  Moreover, given any $p>1$ one can find $C>0$ such that
  \be{14.2}
	\|n(\cdot,t)\|_{L^\infty(\Omega)}
	+ \|c(\cdot,t)\|_{W^{1,p}(\Omega)}
	+ \|u(\cdot,t)\|_{L^\infty(\Omega)} 
	\le C
	\qquad \mbox{for all } t\ge 0.
  \ee
\end{theo}
With regard to the question of global solvability by bounded functions for arbitrary coefficient functions $\phi$ and $f$
and initial data, the problem of identifying a critical decay rate of $S$,
up to evident remaining open topics arising when e.g.~in (\ref{proto}) we precisely have equality,
thereby seems comprehensively solved in the spatially three-dimensional case.
In comparison to this, the picture seems much less complete in neighboring families of systems	
in which chemotactic cross-diffusion interacts with either alternative or further mechanisms.
For instance, logistic-type growth restrictions, as modeled by additional summands of the
form $\rho n - \mu n^2$ in the respective equation for $n$, have recently been shown to prevent blow-up in
corresponding Keller-Segel-fluid variants of (\ref{0}) if either $N=2$ and $\mu>0$ is arbitrary, even in the case when 
the fluid flow is governed by the full Navier-Stokes equations (\cite{espejo_suzuki}, \cite{taowin_ZAMP2016}), 
or $N=3$ and $\mu>0$ is suitably large (\cite{taowin_ZAMP2015}). 
This generalizes previously known facts for the corresponding fluid-free Keller-Segel-growth system (\cite{OTYM},
\cite{win_CPDE2010}), but due to the lack of any complementary result on blow-up e.g.~for $N=3$ and small $\mu>0$,
this only partially clarifies how far the potential to enforce explosions is influenced by fluid interaction 
in such circumstances.
Similar observations concern related chemotaxis(-fluid) systems
accounting for consumption, rather than production, of the chemical signal by the cells,
in the most prototypical form requiring a replacement of the reaction term $-c+n$ with $-nc$ in the equation 
determining the evolution of $c$.
Models of this form have been studied quite thoroughly in the literature, both with diffusion and cross-diffusion
of the form in (\ref{0}) (\cite{DLM}, \cite{win_CPDE2012}, \cite{chae_kang_lee_CPDE}, \cite{cao_lankeit}), 
and also with focus on blow-up-inhibiting effects of 
either nonlinear variants of cross-diffusion rates as in (\ref{KSM}) 
(\cite{wang_cao}, \cite{wangyulan_CAMWA}),		
or of porous medium-type diffusion (\cite{DiFLM}, \cite{duan_xiang},
\cite{taowin_ANIHPC}, 		
\cite{win_ct_fluid_nonlin}, \cite{yilong_wang_xie}).
In fact, various sets of conditions could be identified as sufficient for global solvability
in such systems within classes of bounded functions (\cite{yilong_wang_xie},	
\cite{win_ct_fluid_nonlin}, \cite{win_ARMA}, \cite{taowin_DCDSA}),
but due to missing examples of blow-up it seems widely unclear yet how far they are necessary therefor in the respective
setting.\abs
We remark that as a by-product, Theorem \ref{theo14} also asserts global existence of bounded solutions to the
corresponding Neumann initial-boundary value problem for the two-component chemotaxis-transport system 
\bas
	\left\{ \begin{array}{lcll}
    	n_t + u\cdot\nabla n &=& \Delta n  - \nabla \cdot \Big(nS(n)\nabla c\Big), 
	\qquad & x\in\Omega, \ t>0,\\[1mm]
    	c_t + u\cdot\nabla c &=& \Delta c-c+n, 
	\qquad & x\in\Omega, \ t>0, 
	\end{array} \right.
\eas
with any prescribed sufficiently smooth and bounded solenoidal fluid field $u$;
in fact, this can readily be verified upon obvious choices of $\phi$ and $f$ in Theorem \ref{theo14}.\abs
{\bf Main ideas.} \quad
In the literature on 
the fluid-free system (\ref{KSM}), proofs for boundedness under the optimal version of (\ref{S}) could be built 
on analyzing functionals of the form
\bas
	y(t):=\io n^p(\cdot,t) + \io |\nabla c(\cdot,t)|^r,
	\qquad t>0,
\eas
for suitably chosen $p>1$ and $r>1$ (\cite{horstmann_win}, \cite{taowin_subcrit}). 
Indeed, it can be seen that in a	
correspondingly obtained ODE for $y$,
by making use of (\ref{S}) it becomes possible to control the respective crucial 
cross-diffusive contribution by means of appropriate interpolation in order to show that $y$ satisfies
an ODI of the form $y'+ay \le b$ with some $a>0$ and $b>0$.
However, besides on mass conservation any such interpolation procedure appears to rely on uniform boundedness
of $c$ with respect to the norm in $L^q(\Omega)$ for $q$ close to the largest value $\frac{N}{N-2}$ that can be expected 
for such a property in the heat equation $c_t=\Delta c - c + h$ in $\Omega\times (0,T)$ with $h$ only known to belong
to $L^\infty((0,T);L^1(\Omega))$.\abs
Now in presence of an additional fluid interaction of the form in (\ref{0}), 
it seems unclear whether this is sufficient to warrant that the latter basic
integrability property of the signal $c$ remains to be valid in the entire optimal range $1\le q<3=\frac{N}{N-2}$;
accordingly, pursuing strategies in the flavor of the above needs to cope with weaker a priori information
on $c$ which eventually requires stronger assumptions, 
such as e.g.~in \cite{wang_xiang_JDE2016}, where bounds for $c$ in $L^\infty((0,T);L^2(\Omega))$, yet
available in the whole regime $\alpha>\frac{1}{3}$, are used to finally derive boundedness under the suboptimal 
condition $\alpha>\frac{1}{2}$.\abs
A major technical challenge will thus consist in developing an alternative approach capable of deriving
boundedness of solutions in the optimal range of $\alpha$ but relying on basic regularity information on
$n, c$ and $u$ not substantially going beyond that mentioned above.
In the present work this will be achieved by a series of arguments which at their core are based on an analysis
of the simple functional $z(t):=\io n^p(\cdot,t)$, $t>0$,
for suitably large $p>1$. 
In order to appropriately estimate the respective cross-diffusive summand arising in an associated ODE for $z$
(cf.~(\ref{8.3})),
unlike in most previous related works we shall make essential use 
of maximal Sobolev regularity properties of the heat and the Stokes evolution equations to derive bounds for
the divergence $\Delta c$ of the cross-diffusive gradient 
which immediately arises herein (Lemma \ref{lem4} and Lemma \ref{lem6}), 
and the velocity $u$ to which the regularity of the latter is linked (Lemma \ref{lem5}). 
These estimates will be formulated in terms of the quantities given by 
\bas
	\I(T):=\sup_{t\in [\tau,T-\tau]} \int_t^{t+\tau} \io |\nabla n^\frac{p}{2}|^2
\eas
for suitable $\tau\in (0,1]$ and within suitable ranges of $T>2\tau$,
and a crucial observation will reveal by means of appropriate interpolation arguments (Lemma \ref{lem2},
Lemma \ref{lem3} and Lemma \ref{lem7})
that when merely $\alpha>\frac{1}{3}$,
for sufficiently large $p>1$ these quantities will satisfy inequalities of the form
$\I(T) \le a\I^\gamma(T)+b$ with some $a>0,b>0$ and $\gamma\in (0,1)$ conveniently independent of $T$ (Lemma \ref{lem8}).
The boundedness properties of $\I(T)$ thereby implied will afterwards entail estimates for $n$ 
with respect to the norm in $L^p(\Omega)$ for arbitrarily large $p>1$ (Lemma \ref{lem9}) and thus, through subsequent
applications of basically well-established methods, yield estimates sufficient for the derivation of 
Theorem \ref{theo14} (Section \ref{sect_theo14}).
We emphasize that during our interpolation procedures we shall only rely on an easily obtained weak
a priori boundedness information on 
$(n,c,u)$ in the spaces $L^1(\Omega)\times L^1(\Omega) \times L^p(\Omega;\R^3)$ for arbitrary $p\in (1,3)$
(see Lemma \ref{lem_basic} and Lemma \ref{lem02}).
\mysection{Preliminaries}
\subsection{Local existence and basic solution properties}
Let us first state a basic result on local existence 
and extensibilty that can be achieved by means of arguments well-known
in the theory of chemotaxis and chemotaxis-fluid systems (\cite{win_CPDE2012}, \cite{horstmann_win}, \cite{amann}).
\begin{lem}\label{lem_loc}
  Let $\phi\in C^2(\bom)$, $f\in C^1(\bom\times [0,\infty);\R^3)$ and $S\in C^2([0,\infty))$, 
  and suppose that $n_0, c_0$ and $u_0$ comply with (\ref{init}).
  Then there exist $\tm\in (0,\infty]$ and a uniquely determined quadruple $(n,c,u,P)$ of functions 
  \be{l1}
	\left\{ \begin{array}{l}
	n \in C^0(\bom\times [0,\tm)) \cap C^{2,1}(\bom\times (0,\tm)), \\
	c \in \bigcap_{p>3} C^0([0,\tm);W^{1,p}(\Omega)) \cap C^{2,1}(\bom\times (0,\tm)), \\
	u \in C^0([0,\tm);D(A^\beta)) \cap C^{2,1}(\bom\times [0,\tm);\R^3), \\
	P \in C^{1,0}(\bom\times (0,\tm)),
	\end{array} \right.
  \ee
  which are such that
  $n \ge 0$ and $c \ge 0$ in $\Omega\times (0,\tm)$, that $(n,c,u,P)$ solves (\ref{0}) in the classical sense
  in $\Omega\times (0,\tm)$, and that
  \be{ext}
	\mbox{if $\tm<\infty$ \quad then \quad }
	\limsup_{t\nearrow\tm} \Big( \|n(\cdot,t)\|_{L^\infty(\Omega)}
	+ \|c(\cdot,t)\|_{W^{1,p}(\Omega)} + \|A^\beta u(\cdot,t)\|_{L^2(\Omega)} \Big) = \infty
	\quad \mbox{for all } p>3.
  \ee
\end{lem}
The first two solution components can easily be seen to belong to $L^\infty((0,\tm);L^1(\Omega))$:
\begin{lem}\label{lem_basic}
  Under the assumptions of Lemma \ref{lem_loc}, the solution of (\ref{0}) satisfies
  \be{mass}
	\io n(\cdot,t)=\io n_0 
	\qquad \mbox{for all } t \in [0,\tm)
  \ee
  and
  \be{mass_c}
	\io c(\cdot,t) \le \max \bigg\{ \io n_0 \, , \, \io c_0 \bigg\}
	\qquad \mbox{for all } t \in [0,\tm).
  \ee

\end{lem}
\proof 
  We firstly obtain (\ref{mass}) as an immediate consequence of the fact that $\frac{d}{dt} \io n=0$ for all
  $t\in (0,\tm)$ by (\ref{0}). Thereafter, noting that thus $\frac{d}{dt} \io c=- \io c + \io n=-\io c + \io n_0$ 
  for all $t\in (0,\tm)$, 
  we may invoke an ODE comparison argument to readily verify (\ref{mass_c}).
\qed
Under the boundedness assumption on $f$ from Theorem \ref{theo14}, due to (\ref{mass}) 
also the fluid velocity enjoys a basic boundedness property. 
As precedent derivations of similar features in related systems apparently only address contexts without external
source terms (see e.g.~\cite[Lemma 2.5]{wang_xiang_JDE2016}), let us include a short proof 
of this essentially well-known fact here for completeness.
\begin{lem}\label{lem02}
  If, beyond the assumptions of Lemma \ref{lem_loc}, $f$ is bounded in $\Omega\times (0,\infty)$, then
  for each $p\in (1,3)$ there exists $C>0$ such that
  \be{02.1}
	\|u(\cdot,t)\|_{L^p(\Omega)} \le C
	\qquad \mbox{for all } t\in [0,\tm).
  \ee
\end{lem}
\proof
  Since $p<3$ and hence $\frac{3}{2}-\frac{3}{2p}<1$, it is possible to fix $\gamma\in (0,1)$ such that
  $\gamma>\frac{3}{2}-\frac{3}{2p}$, which by a known embedding property (\cite[Lemma 3.3]{win_CVPDE}) ensures
  the existence of $C_1>0$ such that
  \bas
	\|A^{-\gamma}\proj \varphi\|_{L^p(\Omega)}
	\le C_1 \|\varphi\|_{L^1(\Omega)}
	\qquad \mbox{for all } \varphi\in C^1(\bom;\R^3).
  \eas
  According to well-known smoothing properties of the Stokes semigroup (\cite{sohr}, \cite{giga1986}), on the basis
  of a variation-of-constants 
  representation of $u$ we thus infer that with some $C_2>0$ and $\lambda_1>0$ we have
  \bas
	\|u(\cdot,t)\|_{L^p(\Omega)}
	&=& \bigg\| e^{-tA}u_0 + \int_0^t A^\gamma e^{-(t-s)A} A^{-\gamma} 
		\proj \Big[ n(\cdot,s)\nabla \phi + f(\cdot,s)\Big] ds \bigg\|_{L^p(\Omega)} ds \\
	&\le& C_2 \|u_0\|_{L^p(\Omega)}
	+ C_2 \int_0^t (t-s)^{-\gamma} e^{-\lambda_1(t-s)} 
		\Big\| A^{-\gamma} \proj\Big[n(\cdot,s)\nabla \phi + f(\cdot,s)\Big] \Big\|_{L^1(\Omega)} ds \\
	&\le& C_2\|u_0\|_{L^p(\Omega)}
	+ C_1 C_2 \int_0^t (t-s)^{-\gamma} e^{-\lambda_1(t-s)} 
		\big\| n(\cdot,s)\nabla \phi + f(\cdot,s)\big\|_{L^1(\Omega)}
  \eas
  for all $t\in [0,\tm)$. Since using (\ref{mass}) we obtain that
  \bas
	\big\| n(\cdot,s)\nabla \phi + f(\cdot,s)\big\|_{L^1(\Omega)}
	\le \|\nabla\phi\|_{L^\infty(\Omega)} \io n_0 +|\Omega| \cdot \|f\|_{L^\infty(\Omega\times (0,\infty))}
	\qquad \mbox{for all } s\in (0,\tm),
  \eas
  and since the requirement that $\gamma<1$ implies that
  $\int_0^t (t-s)^{-\gamma} e^{-\lambda_1(t-s)} ds \le \int_0^\infty \sigma^{-\gamma} e^{-\lambda_1 \sigma} d\sigma
  <\infty$
  for all $t\ge 0$, this immediately yields (\ref{02.1}).
\qed
\subsection{An ODE lemma}
For later use in Lemma \ref{lem3} and Lemma \ref{lem8}, let us provide an elementary statement on upper estimates
in superlinearly dampened ordinary differential inequalities involving forcing terms only known to be bounded in average.
\begin{lem}\label{lem1}
  Let $\ts\in \R$, $T>\ts$ and $\tau\in (0,T-\ts)$, and suppose that $y\in C^0([\ts,T)) \cap C^1((\ts,T))$, 
  $g\in L^1((\ts,T))$ and $h\in L^1((\ts,T))$ are nonnegative and such that
  \be{1.1}
	y'(t) + ay^\gamma(t) + g(t) \le h(t)
	\qquad \mbox{for all } t\in (\ts,T)
  \ee
  and
  \be{1.2}
	\int_t^{t+\tau} h(s) ds \le b
	\qquad \mbox{for all } t\in [\ts,T-\tau]
  \ee
  with some $a>0, b>0$ and $\gamma>1$. Then
  \be{1.3}
	y(t) \le b+C
	\qquad \mbox{for all } t\in [\ts,T]
  \ee
  and
  \be{1.4}
	\int_t^{t+\tau} g(s) ds \le 2b+C
	\qquad \mbox{for all } t\in [\ts,T-\tau],
  \ee
  where
  \be{1.5}
	C:=\max \bigg\{ y(\ts) \, , \, [(\gamma-1) a\tau]^{-\frac{1}{\gamma-1}} \bigg\}.
  \ee
\end{lem}
\proof
  Abbreviating $C_1:=[(\gamma-1)a]^{-\frac{1}{\gamma-1}}$ and without loss of generality assuming that $\ts=0$, 
  we first claim that then for any choice of $t_0\in [0,T)$ we 
  have
  \be{1.6}
	y(t) \le \oy(t):=C_1(t-t_0)^{-\frac{1}{\gamma-1}}
	+ \int_{t_0}^t h(s) ds
	\qquad \mbox{for all } t\in (t_0,T].
  \ee
  To verify this, we observe that
  \bas
	\oy'(t) + a\oy^\gamma(t) - h(t)
	&=& -\frac{C_1}{\gamma-1} (t-t_0)^{-\frac{\gamma}{\gamma-1}} +h(t) 
	+ a\cdot \bigg\{ C_1(t-t_0)^{-\frac{1}{\gamma-1}} + \int_{t_0}^t h(s)ds \bigg\}^\gamma - h(t) \\
	&\ge& - \frac{C_1}{\gamma-1} (t-t_0)^{-\frac{\gamma}{\gamma-1}} 
	+ a\cdot \bigg\{ C_1 (t-t_0)^{-\frac{1}{\gamma-1}} \bigg\}^\gamma \\[2mm]
	&=& 0
	\qquad \mbox{for all } t\in (t_0,T),
  \eas
  because $\frac{C_1}{\gamma-1}=aC_1^\gamma$ according to our definition of $C_1$.
  Since $y$ is bounded and $\oy(t)\nearrow +\infty$ as $t\searrow t_0$, an ODE comparison argument on
  $[t_0+\delta,T]$ with suitably small $\delta\in (0,T-t_0)$ therefore yields (\ref{1.6}).\\
  Now for $t\ge \tau$, we may therein choose $t_0:=t-\tau$ to see that in view of (\ref{1.2}) and (\ref{1.5}),
  \bas
	y(t) \le C_1 \tau^{-\frac{1}{\gamma-1}}
	+ \int_{t-\tau}^t h(s) ds
	\le C_1 \tau^{-\frac{1}{\gamma-1}} + b
	\le C+b
	\qquad \mbox{for all } t\in [\tau,T],
  \eas
  whereas for smaller $t$ we simply neglect two nonnegative summands on the left of (\ref{1.1}) to find upon integration
  that again due to (\ref{1.2}),
  \bas
	y(t) \le y(0) + \int_0^t h(s) ds 
	\le C + \int_0^\tau h(s) ds 
	\le C+b 
	\qquad \mbox{for all } t\in [0,\tau),
  \eas
  because $h$ is nonnegative.\\
  Having thereby established (\ref{1.3}), by means of another integration in (\ref{1.1}) we finally obtain that
  \bas
	\int_t^{t+\tau} g(s)ds
	&\le& y(t)-y(t+\tau)+ \int_t^{t+\tau} h(s)ds \\
	&\le& (b+C) + b
	\qquad \mbox{for all } t\in [0,T-\tau],
  \eas
  and that thus also (\ref{1.4}) is valid.
\qed 
\mysection{A space-time regularity property of $n$ implied by bounds for $\nabla n^\frac{p}{2}$}
In order to simplify notation,
throughout the remaining analysis we assume unless otherwise stated that $\phi$, $f$, $S$ and $(n_0,c_0,u_0)$
are such that the hypotheses of Lemma \ref{lem_loc} are satisfied, that moreover $f$ is bounded, and that
(\ref{S}) holds with some $\ks>0$ and $\alpha>0$. We then let $(n,c,u,P)$ and $\tm\in (0,\infty]$
be as provided by Lemma \ref{lem_loc}, and set
\be{tau}
	\tau:=\min\Big\{1 \, , \, \frac{1}{4}\tm\Big\}.
\ee
Now in the major part of our subsequent reasoning, a crucial role will be played by the quantities defined by
\be{I}
	\I(T):=\sup_{t\in [\tau,T-\tau]} \int_t^{t+\tau} \io |\nabla n^\frac{p}{2}|^2
	\qquad \mbox{for } T\in (2\tau,\tm)
	\quad \mbox{and }  p>1,
\ee
which contain the dissipated quantity appearing in a standard $L^p$ testing procedure when applied to the first
equation in (\ref{0}). 
Our arguments to control the cross-diffusive contributions therein will be prepared by a series of bounds
for $n$, $c$ and $u$ in terms of $\I$, with the final ambition to estimate $\I(T)$ by, essentially, a sublinear
power thereof (Lemma \ref{lem8}).\abs
Our first step in this direction, based on a simple interpolation argument involving (\ref{mass}),
will frequently be applied in the following lemmata.
\begin{lem}\label{lem2}
  Let $p>1$, and suppose that $\kappa>1$ and $\lambda>0$ are such that
  \be{2.1}
	\kappa\le 3p
  \ee
  and
  \be{2.2}
	\frac{3(\kappa-1)}{(3p-1)\kappa} \cdot \lambda <1.
  \ee
  Then there exists $C>0$ such that for all $T\in (2\tau,\tm)$ we have
  \be{2.3}
	\int_t^{t+\tau} \|n(\cdot,s)\|_{L^\kappa(\Omega)}^\lambda ds
	\le C + C \I^\frac{3(\kappa-1)\lambda}{(3p-1)\kappa}(T)
	\qquad \mbox{for all } t\in [0,T-\tau].
  \ee
\end{lem}
\proof
  Using that $\kappa\ge 1$ and $\kappa\le 3p$, we invoke the Gagliardo-Nirenberg inequality to fix $C_1>0$ such that
  \bas
	\|\varphi\|_{L^\frac{2\kappa}{p}(\Omega)}^\frac{2\lambda}{p}
	\le C_1\|\nabla\varphi\|_{L^2(\Omega)}^\frac{6(\kappa-1)\lambda}{(3p-1)\kappa}
		\|\varphi\|_{L^\frac{2}{p}(\Omega)}^\frac{2(3p-\kappa)\lambda}{p(3p-1)\kappa}
	+ C_1\|\varphi\|_{L^\frac{2}{p}(\Omega)}^\frac{2\lambda}{p}
	\qquad \mbox{for all } \varphi\in W^{1,2}(\Omega),
  \eas
  and apply this to $\varphi:= n^\frac{p}{2}(\cdot,s)$ for $s\in (\tau,\tm)$ to see upon a time integration that
  since $\|n^\frac{p}{2}(\cdot,s)\|_{L^\frac{2}{p}(\Omega)}^\frac{2}{p} = C_2:=\io n_0$
  for all $s \in (\tau,\tm)$ by (\ref{mass}),
  \bea{2.4}
	\int_t^{t+\tau} \|n(\cdot,s)\|_{L^\kappa(\Omega)}^\lambda ds
	&=& \int_t^{t+\tau} \|n^\frac{p}{2}(\cdot,s)\|_{L^\frac{2\kappa}{p}(\Omega)}^\frac{2\lambda}{p} ds \nn\\
	&\le& C_3 \int_t^{t+\tau} 
		\|\nabla n^\frac{p}{2}(\cdot,s)\|_{L^2(\Omega)}^\frac{6(\kappa-1)\lambda}{(3p-1)\kappa} ds
	+ C_4 \tau
	\qquad \mbox{for all } t\in [\tau,\tm-\tau)
  \eea
  with $C_3:=C_1 C_2^\frac{(3p-\kappa)\lambda}{(3p-1)\kappa}$ and $C_4:= C_1 C_2^\lambda$.
  Here thanks to the fact that $\frac{6(\kappa-1)\lambda}{(3p-1)\kappa}<2$ according to (\ref{2.2}), we may employ
  the H\"older inequality to obtain that
  \bas
	\int_t^{t+\tau} \|\nabla n^\frac{p}{2}(\cdot,s)\|_{L^2(\Omega)}^\frac{6(\kappa+1)\lambda}{(3p-1)\kappa} ds
	\le \bigg\{ \int_t^{t+\tau} \|\nabla n^\frac{p}{2}(\cdot,s)\|_{L^2(\Omega)}^2 ds 
		\bigg\}^\frac{3(\kappa-1)\lambda}{(3p-1)\kappa} \cdot \tau^{1-\frac{3(\kappa-1)\lambda}{(3p-1)\kappa}}
	\qquad \mbox{for all } t\in [\tau,\tm-\tau).
  \eas
  Using that $\tau\le 1$, given $T\in (2\tau,\tm)$ from (\ref{2.4}) we thus conclude that due to the definition of $\I$
  we have
  \bas
	\int_t^{t+\tau} \|n(\cdot,s)\|_{L^\kappa(\Omega)}^\lambda ds
	\le C_3 \I^\frac{3(\kappa-1)\lambda}{(3p-1)\kappa}(T)
	+ C_4
	\qquad \mbox{for all } t\in [\tau,T-\tau],
  \eas
  which implies (\ref{2.3}) due to the fact that $n$ is bounded in $\Omega\times [0,\tau)$ by Lemma \ref{lem_loc}.
\qed
\mysection{An $L^q$ bound for $c$ in terms of $\I(T)$}
A first application of Lemma \ref{lem2} 
yields the following $L^q$ estimate for $c$ in dependence on $\I(T)$, provided that $p$ is suitably large relative to $q$.
Our derivation thereof is based on an $L^q$ testing procedure for the second equation in (\ref{0}) and thus,
due to the solenoidality of the velocity field, does not rely on any explicit bound on $u$.
\begin{lem}\label{lem3}
  Let $p>1$ and $q>1$ be such that
  \be{3.1}
	q<3p.
  \ee
  Then there exists $C>0$ such that for all $T\in (2\tau,\tm)$,
  \be{3.2}
	\io c^q(\cdot,t) \le C + C \I^\frac{q-1}{3p-1}(T)
	\qquad \mbox{for all } t\in [0,T].
  \ee
\end{lem}
\proof
  We use $c^{q-1}$ as a test function for the second equation in (\ref{0}) and note that $\nabla\cdot u\equiv 0$
  to see by means of the H\"older inequality that
  \bea{3.3}
	\frac{1}{q} \frac{d}{dt} \io c^q + (q-1) \io c^{q-2} |\nabla c|^2 + \io c^q
	&=& \io nc^{q-1} = \io n \cdot (c^\frac{q}{2})^\frac{2(q-1)}{q} \nn\\
	&\le& \|n\|_{L^\frac{3q}{2q+1}(\Omega)} \|c^\frac{q}{2}\|_{L^6(\Omega)}^\frac{2(q-1)}{q}
	\qquad \mbox{for all } t\in (0,\tm).
  \eea
  Here employing the three-dimensional Sobolev inequality followed by Young's inequality we can find $C_1>0$ and $C_2>0$
  such that
  \bas	\|n\|_{L^\frac{3q}{2q+1}(\Omega)} \|c^\frac{q}{2}\|_{L^6(\Omega)}^\frac{2(q-1)}{q}
	&\le& C_1\|n\|_{L^\frac{3q}{2q-1}(\Omega)}
	\cdot \bigg\{ \|\nabla c^\frac{q}{2}\|_{L^2(\Omega)}^2 + \|c^\frac{q}{2}\|_{L^2(\Omega)}^2 \bigg\}^\frac{q-1}{q}\\
	&\le& \frac{2(q-1)}{q^2} 
	\cdot \bigg\{ \|\nabla c^\frac{q}{2}\|_{L^2(\Omega)}^2 + \|c^\frac{q}{2}\|_{L^2(\Omega)}^2 \bigg\}
	+ C_2\|n\|_{L^\frac{3q}{2q+1}(\Omega)}^q \\
	&\le& \frac{2(q-1)}{q^2} \io |\nabla c^\frac{q}{2}|^2 + \io c^q + C_2\|n\|_{L^\frac{3q}{2q+1}(\Omega)}^q
	\qquad \mbox{for all } t\in (0,\tm),
  \eas
  because $\frac{2(q-1)}{q^2} \le 1$.
  Since $(q-1) \io c^{q-2}|\nabla c|^2=\frac{4(q-1)}{q^2} \io |\nabla c^\frac{q}{2}|^2$ for all $t\in (0,\tm)$, from
  (\ref{3.3}) we thus infer that
  \be{3.4}
	\frac{d}{dt} \io c^q + C_3 \io |\nabla c^\frac{q}{2}|^2
	\le C_2 q\|n\|_{L^\frac{3q}{2q+1}(\Omega)}^q
	\qquad \mbox{for all } t\in (0,\tm)
  \ee
  with $C_3:=\frac{2(q-1)}{q}>0$. Now combining the Gagliardo-Nirenberg inequality with the fact that
  $\|c^\frac{q}{2}(\cdot,t)\|_{L^\frac{2}{q}(\Omega)}^\frac{2}{q} \le \max\{ \io n_0, \io c_0\}$ by (\ref{mass_c}),
  we can furthermore find $C_4>0$ and $C_5>0$ such that
  \bas
	\bigg\{ \io c^q \bigg\}^\frac{3q-1}{3(q-1)}
	&=& \|c^\frac{q}{2}\|_{L^2(\Omega)}^\frac{2(3q-1)}{3(q-1)} \\
	&\le& C_4\|\nabla c^\frac{q}{2}\|_{L^2(\Omega)}^2 \|c^\frac{q}{2}\|_{L^\frac{2}{q}(\Omega)}^\frac{4}{3(q-1)}
	+ C_4\|c^\frac{q}{2}\|_{L^\frac{2}{q}(\Omega)}^\frac{2(3q-1)}{3(q-1)} \\
	&\le& C_5 \|\nabla c^\frac{q}{2}\|_{L^2(\Omega)}^2 + C_5
	\qquad \mbox{for all } t\in (0,\tm)
  \eas
  and hence
  \bas
	\io |\nabla c^\frac{q}{2}|^2
	\ge \frac{1}{C_5} \cdot \bigg\{ \io c^q \bigg\}^\frac{3q-1}{3(q-1)} -1
	\qquad \mbox{for all } t\in (0,\tm).
  \eas
  Consequently, (\ref{3.4}) can be turned into the inequality
  \bas
	\frac{d}{dt} \io c^q + \frac{C_3}{C_5} \cdot \bigg\{ \io c^q \bigg\}^\frac{3q-1}{3(q-1)}
	\le C_2 q\|n\|_{L^\frac{3q}{2q+1}(\Omega)}^q + C_3
	\qquad \mbox{for all } t\in (0,\tm),
  \eas
  which by Lemma \ref{lem1} implies that whenever $T\in (2\tau,\tm)$,
  \bea{3.5}
	\io c^q(\cdot,t)
	&\le& \sup_{s\in [\tau,T-\tau]} \int_s^{s+\tau} \Big\{ 
	C_2 q\|n(\cdot,\sigma)\|_{L^\frac{3q}{2q+1}(\Omega)}^q + C_3 \Big\} d\sigma +C_6 \nn\\
	&=& C_2 q \sup_{s\in [\tau,T-\tau]} \int_s^{s+\tau} \|n(\cdot,\sigma)\|_{L^\frac{3q}{2q+1}(\Omega)}^q d\sigma
	+ C_3\tau+C_6
	\qquad \mbox{for all } t\in [\tau,T]
  \eea
  with $C_6:=\max \Big\{ \io c_0^q , [\frac{2}{3(q-1)} \cdot \frac{C_3}{C_5} \cdot \tau]^{-\frac{3(q-1)}{2}}\Big\}$.\\
  In order to further estimate the right-hand side in (\ref{3.5}) on the basis of Lemma \ref{lem2}, we observe that
  $\kappa:=\frac{3q}{2q+1}$ and $\lambda:=q$ satisfy
  \bas
	1<\frac{2q+q}{2q+1} =\kappa<\frac{3q}{2q}=\frac{3}{2}<3p
  \eas
  due to our assumptions that $q>1$ and $p>1$, and that our additional requirement (\ref{3.1}) ensures that
  \bas
	\frac{3(\kappa-1)}{(3p-1)\kappa}\cdot\lambda 
	= \frac{q-1}{3p-1}<1.
  \eas
  Therefore, Lemma \ref{lem2} indeed becomes applicable so as to yield $C_7>0$ such that
  \bas
	\int_s^{s+\tau} \|n(\cdot,\sigma\|_{L^\frac{3q}{2q+1}(\Omega)}^q d\sigma
	\le C_7 + C_7 \I^\frac{q-1}{3p-1}(T)
	\qquad \mbox{for all } s\in [\tau,T-\tau],
  \eas
  whereupon (\ref{3.2}) results from (\ref{3.5}) and the boundedness of $c$ in $\Omega\times [0,\tau)$ entailed
  by Lemma \ref{lem_loc}.
\qed
\mysection{Estimates for $\Delta c$ in terms of $\I(T)$ via maximal Sobolev regularity}
Approaching the core of our analysis, our next goal consists in controlling the
cross-diffusive gradient in (\ref{0}) by quantities containing suitably small powers of $\I(T)$ under appropriate
further assumptions on $\alpha$ and the yet free parameter $p$.
Here our first result will relate a second-order Sobolev norm of $c$ to regularity properties
of the three quantities $n,c$ and $u$ making up the inhomogeneity $h:=n-u\cdot\nabla c$ in the heat
equation $c_t=\Delta c - c + h$.
This will be achieved through an argument based on a maximal Sobolev regularity feature of the latter, along with
a suitable temporal regularization procedure which we prepare by
fixing a nondecreasing
$\zeta_0\in C^\infty(\R)$ such that $\zeta_0\equiv 0$ in $(-\infty,-\tau]$ and $\zeta_0\equiv 1$ in $[0,\infty)$,
and defining a family of functions $(\zeta^{(t_0)})_{t_0\in\R}$ by letting
\be{zeta}
	\zeta^{(t_0)}(t):=\zeta_0(t-t_0)
	\qquad \mbox{for $t_0\in\R$ and } t\in\R.
\ee
Our first step toward estimating $\Delta c$ will now consist in the following inequality.
\begin{lem}\label{lem4}
  Let $q>\frac{3}{2}$ and $r>1$. Then there exists $C>0$ such that
  \bea{4.1}
	\int_t^{t+\tau} \|c(\cdot,s)\|_{W^{2,\frac{3}{2}}(\Omega)}^r ds
	&\le& C + C \int_{t-\tau}^{t+\tau} \|n(\cdot,s)\|_{L^\frac{3}{2}(\Omega)}^r ds \nn\\
	& & + C\cdot \bigg\{ \sup_{s\in [0,t+\tau]} \|c(\cdot,s)\|_{L^q(\Omega)} \bigg\}^r
	\cdot \int_{t-\tau}^{t+\tau} \|u(\cdot,s)\|_{L^\frac{6q}{2q-3}(\Omega)}^{2r} ds \nn\\[2mm]
	& & + C \cdot \sup_{s\in [0,t+\tau]} \|c(\cdot,s)\|_{L^\frac{3}{2}(\Omega)}^r
	\qquad \mbox{for all } t\in [\tau,\tm-\tau).
  \eea
\end{lem}
\proof
  We fix $t_0\in [\tau,\tm-\tau)$, and with $\zeta:=\zeta^{(t_0)}$ as defined in (\ref{zeta}) we let
  \bas
	\hc(x,t):=\zeta(t) c(x,t),
	\qquad x\in\bom, \, t\in [t_0-\tau,\tm).
  \eas
  Then $\hc$ is a solution of
  \bas
	\left\{ \begin{array}{ll}
	\hc_t = \Delta\hc -\hc + \zeta n - \zeta u\cdot\nabla c + \zeta_t c,
	\qquad & x\in\Omega, \ t\in (t_0-\tau,\tm), \\[1mm]
	\frac{\partial\hc}{\partial\nu}=0,
	\qquad & x\in\pO, \ t\in (t_0-\tau,\tm),  \\[1mm]
	\hc(x,t_0-\tau)=0,
	\qquad & x\in\Omega,
	\end{array} \right.
  \eas
  so that known results on maximal Sobolev regularity in the Neumann problem for the heat equation 
  (\cite{giga_sohr}) provide $C_1>0$
  such that
  \bea{4.2}
	\int_{t_0-\tau}^{t_0+\tau} \|\hc(\cdot,s)\|_{W^{2,\frac{3}{2}}(\Omega)}^r ds
	&\le& C_1 \int_{t_0-\tau}^{t_0+\tau}
	\Big\| \zeta(s) n(\cdot,s) - \zeta(s) u(\cdot,s)\cdot\nabla c(\cdot,s)
	+ \zeta_t(s) c(\cdot,s) \Big\|_{L^\frac{3}{2}(\Omega)}^r ds \nn\\
	&\le& C_1 \int_{t_0-\tau}^{t_0+\tau} \|n(\cdot,s)\|_{L^\frac{3}{2}(\Omega)}^r ds
	+ C_1 \int_{t_0-\tau}^{t_0+\tau} 
		\Big\| \zeta(s) u(\cdot,s)\cdot\nabla c(\cdot,s)\Big\|_{L^\frac{3}{2}(\Omega)}^r ds \nn\\
	& & + 2C_1 C_2 \sup_{s\in [0,t_0+\tau]} \|c(\cdot,s)\|_{L^\frac{3}{2}(\Omega)}^r
  \eea
  with $C_2:=\|(\zeta_0)_t\|_{L^\infty((-\tau,\tau))}^r$, because $|\zeta_0|\le 1$ and $\tau\le 1$.
  Moreover, using the H\"older inequality we see that
  \bea{4.5}
	& & \hspace*{-30mm}
	\int_{t_0-\tau}^{t_0+\tau} \Big\|\zeta(s) u(\cdot,s)\cdot\nabla c(\cdot,s)\Big\|_{L^\frac{3}{2}(\Omega)}^r ds 
	\nn\\
	&\le& \int_{t_0-\tau}^{t_0+\tau} \zeta^r(s) \|u(\cdot,s)\|_{L^\frac{6q}{2q-3}(\Omega)}^r
	\|\nabla c(\cdot,s)\|_{L^\frac{6q}{2q+3}(\Omega)}^r ds \nn\\
	&\le& \bigg\{ \int_{t_0-\tau}^{t_0+\tau} \|u(\cdot,s)\|_{L^\frac{6q}{2q-3}(\Omega)}^{2r} ds \bigg\}^\frac{1}{2}
	\cdot \bigg\{ \int_{t_0-\tau}^{t_0+\tau} \zeta^{2r}(s) \|\nabla c(\cdot,s)\|_{L^\frac{6q}{2q+3}(\Omega)}^{2r}
		ds \bigg\}^\frac{1}{2},
  \eea
  where according to the Gagliardo-Nirenberg inequality there exists $C_2>0$ such that
  \bas
	\zeta^{2r}(s) \|\nabla c(\cdot,s)\|_{L^\frac{6q}{2q+3}(\Omega)}^{2r}
	&\le& C_2 \zeta^{2r}(s) \|c(\cdot,s)\|_{W^{2,\frac{3}{2}}(\Omega)}^r \|c(\cdot,s)\|_{L^q(\Omega)}^r \\
	&\le& C_2\|\hc(\cdot,s)\|_{W^{2,\frac{3}{2}}(\Omega)}^r
	\cdot \bigg\{ \sup_{\sigma\in [0,t_0+\tau]} \|c(\cdot,\sigma)\|_{L^q(\Omega)} \bigg\}^r
	\quad \mbox{for all } s\in [t_0-\tau,t_0+\tau],
  \eas
  again due to the fact that $|\zeta_0|\le 1$.
  Upon an application of Young's inequality, (\ref{4.5}) therefore entails that
  \bas
	& & \hspace*{-10mm}
	C_1 \int_{t_0-\tau}^{t_0+\tau} \Big\|\zeta(s) u(\cdot,s)\cdot\nabla c(\cdot,s)\Big\|_{L^\frac{3}{2}(\Omega)}^r ds
	\nn\\
	&\le& C_1\sqrt{C_2} \cdot \bigg\{ \sup_{s\in [0,t_0+\tau]} \|c(\cdot,s)\|_{L^q(\Omega)} \bigg\}^\frac{r}{2}
	\cdot \bigg\{ \int_{t_0-\tau}^{t_0+\tau} \|u(\cdot,s)\|_{L^\frac{6q}{2q-3}(\Omega)}^{2r} ds \bigg\}^\frac{1}{2}
	\cdot \bigg\{ \int_{t_0-\tau}^{t_0+\tau} \|\hc(\cdot,s)\|_{W^{2,\frac{3}{2}}(\Omega)}^r ds \bigg\}^\frac{1}{2}
	\\
	&\le& \frac{1}{2} \int_{t_0-\tau}^{t_0+\tau} \|\hc(\cdot,s)\|_{W^{2,\frac{3}{2}}(\Omega)}^r ds
	+ \frac{C_1^2 C_2}{2} \cdot \bigg\{ \sup_{s\in [0,t_0+\tau]} \|c(\cdot,s)\|_{L^q(\Omega)} \bigg\}^r
	\cdot\int_{t_0-\tau}^{t_0+\tau} \|u(\cdot,s)\|_{L^\frac{6q}{2q-3}(\Omega)}^{2r} ds.
  \eas
  In view of the fact that $\zeta\equiv 1$ on $[t_0,t_0+\tau]$ and that hence
  $\hc\equiv c$ throughout $\Omega\times [t_0,t_0+\tau]$, together with (\ref{4.2}) this establishes (\ref{4.1}).
\qed
The expressions on the right of (\ref{4.1}) containing $n$ and $c$ can be estimated in terms of 
$\I(T)$ by means of Lemma \ref{lem2} and Lemma \ref{lem3}. 
In relating the remaining rightmost integral therein to $\I(T)$ as well, 
we rely on a maximal Sobolev regularity property now of the Stokes evolution system to see that this indeed
is possible when the summability power $r$ in (\ref{4.1}) is suitably small.
\begin{lem}\label{lem5}
  Let $p>1, q>\frac{3}{2}$ and $r>1$ be such that
  \be{5.1}
	r<\frac{(3p-1)q}{3}.
  \ee
  Then one can find $C>0$ with the property that for all $T\in (2\tau,\tm)$,
  \be{5.2}
	\int_t^{t+\tau} \|u(\cdot,s)\|_{L^\frac{6q}{2q-3}(\Omega)}^{2r} ds
	\le C + C \I^\frac{3r}{(3p-1)q}(T)
	\qquad \mbox{for all } t\in [0,T-\tau].
  \ee
\end{lem}
\proof
  As $u$ is bounded in $\Omega\times [0,\tau)$ by Lemma \ref{lem_loc}, we only need to derive the claimed inequality
  in the time interval $[\tau,T-\tau]$ for arbitrary $T\in (2\tau,\tm)$.
  To this end, fixing $t_0\in [\tau,T-\tau]$ we once more take $\zeta:=\zeta^{(t_0)}$ from (\ref{zeta}) and let
  \bas
	\hu(x,t):=\zeta(t)u(x,t),
	\qquad x\in\bom, \ t\in [t_0-\tau,\tm),
  \eas
  so that
  \bas
	\left\{
	\begin{array}{ll}
	\hu_t = \Delta \hu - \zeta \nabla P + \zeta n\nabla\phi + \zeta_t u + \zeta f,
	\qquad & x\in\Omega, \ t\in (t_0-\tau,\tm), \\[1mm]
	\hu=0,
	\qquad & x\in\pO, \ t\in (t_0-\tau,\tm), \\[1mm]
	\hu(x,t_0-\tau)=0,
	\qquad & x\in\Omega.
	\end{array} \right.
  \eas
  A maximal Sobolev regularity property of the Stokes evolution semigroup (\cite{giga_sohr}) thus yields $C_1>0$ 
  such that
  \bea{5.3}
	\hspace*{-5mm}
	\int_{t_0-\tau}^{t_0+\tau} \|\hu(\cdot,s)\|_{W^{2,\frac{2q}{2q-1}}(\Omega)}^{2r} ds
	&\le& C_1 \int_{t_0-\tau}^{t_0+\tau} \Big\|\zeta(s) n(\cdot,s)\nabla \phi
	+ \zeta_t(s) u(\cdot,s) + \zeta(s) f(\cdot,s)\Big\|_{L^\frac{2q}{2q-1}(\Omega)}^{2r} ds \nn\\
	& & \hspace*{-15mm}
	\le \ C_2 \int_{t_0-\tau}^{t_0+\tau} \|n(\cdot,s)\|_{L^\frac{2q}{2q-1}(\Omega)}^{2r} ds
	+ C_3 \int_{t_0-\tau}^{t_0+\tau} \|u(\cdot,s)\|_{L^\frac{2q}{2q-1}(\Omega)}^{2r} ds
	+ C_4
  \eea
  with $C_2:= C_1\|\nabla\phi\|_{L^\infty(\Omega)}$, $C_3:=C_1\|(\zeta_0)_t\|_{L^\infty((-\tau,\tau))}$
  and $C_4:=2C_1|\Omega|^\frac{r(2q-1)}{q} \|f\|_{L^\infty(\Omega\times (0,\infty))}^{2r}$, 
  because $|\zeta_0|\le 1$ and $(t_0+\tau)-(t_0-\tau)=2\tau\le 2$.\\
  To estimate the two integrals on the right-hand side herein, we write
  $\kappa:=\frac{2q}{2q-1}$ and $\lambda:=2r$ and note that since $q>\frac{3}{2}$ we have
  \be{5.4}
	1<\kappa<\frac{2\cdot\frac{3}{2}}{2\cdot\frac{3}{2}-1}=\frac{3}{2}<3<3p,
  \ee
  and that thanks to (\ref{5.1}) we moreover know that
  \be{5.5}
	\frac{3(\kappa-1)}{(3p-1)\kappa} \cdot\lambda = \frac{3r}{(3p-1)q} <1.
  \ee
  From (\ref{5.4}) we particularly see that Lemma \ref{lem02} becomes applicable to show that there exists $C_4>0$
  such that
  \be{5.6}
	\|u(\cdot,s)\|_{L^\frac{2q}{2q-1}(\Omega)} \le C_4
	\qquad \mbox{for all } s\in (0,\tm),
  \ee
  and combining (\ref{5.4}) with (\ref{5.5}) we may invoke Lemma \ref{lem2} to find $C_5>0$ fulfilling
  \bas
	\int_t^{t+\tau} \|n(\cdot,s)\|_{L^\frac{2q}{2q-1}(\Omega)}^{2r} ds
	\le C_5 + C_5 \I^\frac{3r}{(3p-1)q}(T)
	\qquad \mbox{for all } t\in [0,T-\tau],
  \eas
  so that from (\ref{5.3}) and (\ref{5.6}) we thus infer that
  \bas
	\int_{t_0-\tau}^{t_0+\tau} \|\hu(\cdot,s)\|_{W^{2,\frac{2q}{2q-1}}(\Omega)}^{2r} ds
	\le 2C_2 C_5 + 2C_2 C_5 \I^\frac{3r}{(3p-1)q}(T)
	+ 2C_3 C_4^{2r} \tau,
  \eas
  Since $\hu\equiv u$ in $\Omega\times [t_0,t_0+\tau]$ by (\ref{zeta}), and since
  $W^{2,\frac{2q}{2q-1}}(\Omega) \hra L^\frac{6q}{2q-3}(\Omega)$ in the present three-dimensional setting,
  this establishes (\ref{5.2}) in the time interval $[\tau,T-\tau]$, as intended.
\qed 
We can now formulate the main result of this section by combining Lemma \ref{lem4} with Lemma \ref{lem5}, Lemma \ref{lem3}
and Lemma \ref{lem2}, where the latter turns out to be applicable here under a further smallness assumption on $r$.
\begin{lem}\label{lem6}
  Suppose that $p>1, q>\frac{3}{2}$ and $r>1$ are such that (\ref{3.1}) and (\ref{5.1}) hold as well as
  \be{6.2}
	r<3p-1.
  \ee
  Then there exists $C>0$ such that whenever $T\in (2\tau,\tm)$,
  \be{6.3}
	\int_t^{t+\tau} \|\Delta c(\cdot,s)\|_{L^\frac{3}{2}(\Omega)}^r ds
	\le C + C \I^\frac{r}{3p-1}(T)
	+ C\I^\frac{(q+2)r}{(3p-1)q}(T)
	\qquad \mbox{for all } t\in [\tau,T-\tau].
  \ee
\end{lem}
\proof
  Based on our assumptions that $q>\frac{3}{2}$ and that (\ref{3.1}) and (\ref{5.1}) hold, we first employ 
  Lemma \ref{lem3} and Lemma \ref{lem5} to find positive constants $C_1, C_2$ and $C_3$ such that given $T\in (2\tau,\tm)$ we know that
  \be{6.4}
	\|c(\cdot,s)\|_{L^q(\Omega)}^r
	\le C_1 + C_1 \I^\frac{(q-1)r}{(3p-1)q}(T)
	\qquad \mbox{for all } s\in [0,T]
  \ee
  and
  \be{6.44}
	\|c(\cdot,s)\|_{L^\frac{3}{2}(\Omega)}^r
	\le C_2 + C_2 \I^\frac{r}{3(3p-1)}(T)
	\qquad \mbox{for all } s\in [0,T]
  \ee
  as well as
  \be{6.5}
	\int_{t_0}^{t_0+\tau} \|u(\cdot,s)\|_{L^\frac{6q}{2q-3}(\Omega)}^{2r} ds
	\le C_3 + C_3 \I^\frac{3r}{(3p-1)q} (T)
	\qquad \mbox{for all } t_0\in [0,T-\tau].
  \ee
  Moreover, writing $\kappa:=\frac{3}{2}$ and $\lambda:=r$ we obviously have $1<\kappa<3p$, whereas (\ref{6.2}) guarantees
  that
  \bas
	\frac{3(\kappa-1)}{(3p-1)\kappa}\cdot \lambda<1,
  \eas
  so that as a consequence of Lemma \ref{lem2} we can pick $C_4>0$ satisfying
  \be{6.6}
	\int_{t_0}^{t_0+\tau} \|n(\cdot,s)\|_{L^\frac{3}{2}(\Omega)}^r ds
	\le C_4 + C_4\I^\frac{r}{3p-1}(T)
	\qquad \mbox{for all } t_0\in [0,T-\tau].
  \ee
  Now from Lemma \ref{lem4} it follows that there exists $C_5>0$ such that
  \bas
	& & \hspace*{-20mm}
	\int_t^{t+\tau} \|\Delta c(\cdot,s)\|_{L^\frac{3}{2}(\Omega)}^r ds \\
	&\le& C_5
	+ C_5 \int_{t-\tau}^{t+\tau} \|n(\cdot,s)\|_{L^\frac{3}{2}(\Omega)}^r ds \\
	& & + C_5 \cdot \bigg\{ \sup_{s\in [0,t+\tau]} \|c(\cdot,s)\|_{L^q(\Omega)} \bigg\}^r 
	\cdot \int_{t-\tau}^{t+\tau} \|u(\cdot,s)\|_{L^\frac{6q}{2q-3}(\Omega)}^{2r} ds \\
	& & + C_5 \sup_{s\in [0,t+\tau]} \|c(\cdot,s)\|_{L^\frac{3}{2}(\Omega)}
	\qquad \mbox{for all } t\in [\tau,\tm-\tau),
  \eas
  which in light of (\ref{6.4})-(\ref{6.6}) particularly entails that
  \bas
	& & \hspace*{-20mm}
	\int_t^{t+\tau} \|\Delta c(\cdot,s)\|_{L^\frac{3}{2}(\Omega)}^r ds \\
	&\le& C_5 + C_5 \cdot \Big\{ 2C_4 + 2C_4 \I^\frac{r}{3p-1}(T)\Big\} \\
	& & + C_5\cdot \Big\{ C_1 + C_1 \I^\frac{(q-1)r}{(3p-1)q}(T) \Big\}
	\cdot \Big\{ 2C_3 + 2C_3 \I^\frac{3r}{(3p-1)q}(T)\Big\} \\
	& & + C_5 \cdot \Big\{ C_2 + C_2 \I^\frac{r}{3(3p-1)}(T) \Big\}
	\qquad \mbox{for all } t\in [\tau,T-\tau].
  \eas
  As three applications of Young's inequality show that
  \bas
	& & \hspace*{-40mm}
	\Big\{ C_1 + C_1 \I^\frac{(q-1)r}{(3p-1)q}(T) \Big\}
	\cdot \Big\{ 2C_3 + 2C_3 \I^\frac{3r}{(3p-1)q}(T)\Big\} \\
	&=& 2C_1 C_3 \cdot \Big\{ 1 + \I^\frac{(q-1)r}{(3p-1)q}(T)
	+ \I^\frac{3r}{(3p-1)q}(T) + \I^\frac{(q+2)r}{(3p-1)q}(T)\Big\} \\
	&\le& 2C_1 C_3 \cdot \Big\{ 3 + 3\I^\frac{(q+2)r}{(3p-1)q}(T) \Big\}
  \eas
  and
  \bas
	C_2 \I^\frac{r}{3(3p-1)}(T) \le C_2 \I^\frac{r}{3p-1}(T) + C_2,
  \eas
  the derivation of (\ref{6.3}) is complete.
\qed
\mysection{$L^p$ bounds for $n$ by closing the loop}
Now controlling the cross-diffusive action in the announced testing procedure for $n$, to be detailed in
Lemma \ref{lem8}, will amount to appropriately estimating $\io n^{p-\alpha} |\Delta c|$. 
This can be achieved by means of Lemma \ref{lem6} and, again, Lemma \ref{lem2} if the exponent $r$ 
in addition to the assumptions therein satisfies a further condition requiring $r$ not to be too small:
\begin{lem}\label{lem7}
  Let $p>1, q>\frac{3}{2}$ and $r>1$ satisfy (\ref{3.1}), (\ref{5.1}) and (\ref{6.2}) as well as
  \be{7.01}
	p>\alpha+\frac{1}{3}
  \ee
  and
  \be{7.1}
	r>\frac{3p-1}{3\alpha}.
  \ee
  Then one can find $C>0$ such that for each $T\in (2\tau,\tm)$,
  \be{7.2}
	\int_t^{t+\tau} \|n(\cdot,s)\|_{L^{3(p-\alpha)}(\Omega)}^{p-\alpha} 
		\|\Delta c(\cdot,s)\|_{L^\frac{3}{2}(\Omega)} ds
	\le C + C\I^\frac{3(p-\alpha)}{3p-1}(T)
	+ C \I^\frac{3(p-\alpha)q+2}{(3p-1)q}(T)
	\qquad \mbox{for all } t\in [\tau,T-\tau].
  \ee
\end{lem}
\proof
  By the H\"older inequality,
  \bea{7.3}
	& & \hspace*{-18mm}
	\int_t^{t+\tau} \|n(\cdot,s)\|_{L^{3(p-\alpha)}(\Omega)}^{p-\alpha}
	\|\Delta c(\cdot,s)\|_{L^\frac{3}{2}(\Omega)} ds \nn\\
	& & \hspace*{-5mm}
	\le \bigg\{ \int_t^{t+\tau} \|n(\cdot,s)\|_{L^{3(p-\alpha)}(\Omega)}^\frac{(p-\alpha)r}{r-1} ds 
		\bigg\}^{1-\frac{1}{r}}
	\cdot \bigg\{ \int_t^{t+\tau} \|\Delta c(\cdot,s)\|_{L^\frac{3}{2}(\Omega)}^r \bigg\}^\frac{1}{r}
	\ \mbox{for all } t\in [\tau,\tm-\tau).
  \eea
  Here letting $\kappa:=3(p-\alpha)$ we trivially have $\kappa<3p$, while (\ref{7.01}) asserts that $\kappa>1$.
  Furthermore, the hypothesis (\ref{7.1}) guarantees that if we define $\lambda:=\frac{(p-\alpha)r}{r-1}$, then
  \bas
	\frac{3(\kappa-1)}{(3p-1)\kappa}\cdot\lambda
	= \frac{3(p-\alpha)-1}{(3p-1)\cdot (1-\frac{1}{r})}
	< \frac{3(p-\alpha)-1}{(3p-1)\cdot (1-\frac{3\alpha}{3p-1})}
	=1,
  \eas
  whence invoking Lemma \ref{lem2} we can fix $C_1>0$ such that for all $T\in (2\tau,\tm)$,
  \bea{7.4}
	\bigg\{ \int_t^{t+\tau} \|n(\cdot,s)\|_{L^{3(p-\alpha)}(\Omega)}^\frac{(p-\alpha)r}{r-1} ds 
		\bigg\}^{1-\frac{1}{r}}
	&\le& C_1 + C_1\I^{\frac{3(p-\alpha)-1}{(3p-1)\cdot (1-\frac{1}{r})} \cdot (1-\frac{1}{r})}(T) \nn\\
	&=& C_1 + C_1 \I^\frac{3(p-\alpha)-1}{3p-1}(T)
	\qquad \mbox{for all } t\in [\tau,T-\tau].
  \eea
  Next, relying on (\ref{3.1}), (\ref{5.1}) and (\ref{6.2}) we employ Lemma \ref{lem6} to find $C_2>0$ with the property
  that  for any such $T$,
  \bas
	\bigg\{ \int_t^{t+\tau} \|\Delta c(\cdot,s)\|_{L^\frac{3}{2}(\Omega)}^r ds \bigg\}^\frac{1}{r}
	\le C_2 + C_2\I^\frac{1}{3p-1}(T)
	+ C_2 \I^\frac{q+2}{(3p-1)q}(T)
	\qquad \mbox{for all } t\in [\tau,T-\tau].
  \eas
  In conjunction with (\ref{7.3}) and (\ref{7.4}), on three straightforward applications of Young's inequality this
  shows that
  \bas
	& & \hspace*{-20mm}
	\int_t^{t+\tau} \|n(\cdot,s)\|_{L^{3(p-\alpha)}(\Omega)}^{p-\alpha}
	\|\Delta c(\cdot,s)\|_{L^\frac{3}{2}(\Omega)} ds \nn\\
	&\le& C_1 C_2 \cdot \Big\{ 1 + \I^\frac{3(p-\alpha)-1}{3p-1}(T)\Big\} \cdot
	\Big\{ 1+\I^\frac{1}{3p-1}(T)+ \I^\frac{q+2}{(3p-1)q}(T)\Big\} \\
	&=& C_1 C_2\cdot \Big\{ 1 + \I^\frac{1}{3p-1}(T) + \I^\frac{q+2}{(3p-1)q}(T) 
	+ \I^\frac{3(p-\alpha)-1}{3p-1}(T) + \I^\frac{3(p-\alpha)}{3p-1}(T)
	+ \I^\frac{3(p-\alpha)q+2}{(3p-1)q}(T) \Big\} \\
	&\le& C_1 C_2 \cdot \Big\{ 4 + 3\I^\frac{3(p-\alpha)}{3p-1}(T) + 2\I^\frac{3(p-\alpha)q+2}{(3p-1)q}(T)\Big\}
  \eas
  for all $t\in [\tau,T-\tau]$.
\qed
We are now prepared for closing our circle of arguments by an application of Lemma \ref{lem1} to an ODI obtained 
on the basis of the announced $L^p$ testing when combined with Lemma \ref{lem7},
provided that $\alpha$ satisfies the assumption from Theorem \ref{theo14} and the exponent $q$ originating from 
Lemma \ref{lem3} is thereafter fixed appropriately large.
\begin{lem}\label{lem8}
  Suppose that $\alpha>\frac{1}{3}$, and let $p>1, q>\frac{3}{2}$ and $r>1$ be such that (\ref{3.1}), (\ref{5.1}),
  (\ref{6.2}), (\ref{7.01}) and (\ref{7.1}) hold, and such that moreover
  \be{8.1}
	q>\frac{2}{3\alpha-1}.
  \ee
  Then there exists $C>0$ such that
  \be{8.2}
	\io n^p(\cdot,t) \le C
	\qquad \mbox{for all } t\in [0,\tm).
  \ee
\end{lem}
\proof
  We multiply the first equation in (\ref{0}) by $n^{p-1}$ to find using several integrations by parts that
  \bea{8.3}
	\frac{d}{dt} \io n^p
	+ p(p-1)\io n^{p-2} |\nabla n|^2
	&=& p(p-1) \io n^{p-1} S(n) \nabla n\cdot\nabla c \nn\\
	&=& p(p-1) \io \nabla\Psi(n)\cdot\nabla c \nn\\
	&=& -p(p-1) \io \Psi(n)\Delta c
	\qquad \mbox{for all } t\in (0,\tm),
  \eea
  where we have set
  \bas
	\Psi(\xi):=\int_0^\xi \sigma^{p-1} S(\sigma) d\sigma
	\qquad \mbox{for } \xi\ge 0.
  \eas
  Here thanks to (\ref{S}), we can estimate
  \bas
	|\Psi(\xi)| 
	&\le& \ks \int_0^\xi \sigma^{p-1} (\sigma+1)^{-\alpha} d\sigma \\
	&\le& \ks \int_0^\xi \sigma^{p-1-\alpha} d\sigma \\
	&=& \frac{\ks}{p-\alpha} \xi^{p-\alpha}
	\qquad \mbox{for all } \xi\ge 0,
  \eas
  so that by means of the H\"older inequality, on the right-hand side of (\ref{8.3}) we obtain
  \bas
	-p(p-1) \io \Psi(n)\Delta c
	&\le& C_1 \io n^{p-\alpha} |\Delta c| \\
	&\le& C_1\|n\|_{L^{3(p-\alpha)}(\Omega)}^{p-\alpha} \|\Delta c\|_{L^\frac{3}{2}(\Omega)}
	\qquad \mbox{for all } t\in (0,\tm)
  \eas
  with $C_1:=\frac{p(p-1)\ks}{p-\alpha}$.
  Apart from that, using the Gagliardo-Nirenberg inequality together with (\ref{mass}) we see that with some $C_2>0$
  and $C_3>0$ we have
  \bas
	\bigg\{ \io n^p \bigg\}^\frac{3p-1}{3(p-1)}
	=\|n^\frac{p}{2}\|_{L^2(\Omega)}^\frac{2(3p-1)}{3(p-1)}
	\le C_2 \|\nabla n^\frac{p}{2}\|_{L^2(\Omega)}^2 \|n^\frac{p}{2}\|_{L^\frac{2}{p}(\Omega)}^\frac{4}{3(p-1)}
	+ C_2\|n^\frac{p}{2}\|_{L^\frac{2}{p}(\Omega)}^\frac{2(3p-1)}{3(p-1)}
	\le C_3\|\nabla n^\frac{p}{2}\|_{L^2(\Omega)}^2 + C_3
  \eas
  for all $t\in (0,\tm)$, and that abbreviating $C_4:=\frac{2(p-1)}{p}$ we thus can estimate
  \bas
	p(p-1) \io n^{p-2}|\nabla n|^2
	= 2C_4 \io |\nabla n^\frac{p}{2}|^2
	\ge \frac{C_4}{C_3} \cdot \bigg\{ \io n^p \bigg\}^\frac{3p-1}{3(p-1)} - \frac{C_4}{C_3}
	+ C_4 \io |\nabla n^\frac{p}{2}|^2
  \eas
  for all $t\in (0,\tm)$.
  From (\ref{8.3}) we thus infer that
  \bas
	& & y(t):=\io n^p(\cdot,t),
	\quad 
	g(t):=C_4\io |\nabla n^\frac{p}{2}(\cdot,t)|^2
	\quad \mbox{and} \quad \\
	& & h(t):=\frac{C_4}{C_3} + 
	C_1\|n(\cdot,t)\|_{L^{3(p-\alpha)}(\Omega)}^{p-\alpha} \|\Delta c(\cdot,t)\|_{L^\frac{3}{2}(\Omega)},
	\qquad t\in [\tau,\tm),
  \eas
  satisfy
  \bas
	y'(t) + \frac{C_4}{C_3} y^\frac{3p-1}{3(p-1)}(t) + g(t)\le h(t)
	\qquad \mbox{for all } t\in (\tau,\tm),
  \eas
  where due to (\ref{3.1}), (\ref{5.1}), (\ref{6.2}), (\ref{7.01}) and (\ref{7.1}) we may invoke Lemma \ref{lem7}
  to find $C_5>0$ such that for all $T\in (2\tau,\tm)$ we have
  \be{8.4}
	\int_t^{t+\tau} h(s) ds
	\le C_5 + C_5\I^\frac{3(p-\alpha)}{3p-1}(T)
	+ C_5\I^\frac{3(p-\alpha)q+2}{(3p-1)q}(T)
	\qquad \mbox{for all } t\in [\tau,T-\tau].
  \ee
  Therefore, Lemma \ref{lem1} firstly states that if we let
  \bas
	C_6:=\max \bigg\{ \io n^p(\cdot,\tau) \, , \, 
	\Big[\frac{2}{3(p-1)} \cdot \frac{C_4}{C_3} \cdot \tau\Big]^{-\frac{3(p-1)}{2}}\bigg\},
  \eas
  then 
  \bas
	\int_t^{t+\tau} g(s)ds
	\le 2C_5 + 2C_5\I^\frac{3(p-\alpha)}{3p-1}(T)
	+ 2C_5 \I^\frac{3(p-\alpha)q+2}{(3p-1)q}(T) + C_6
	\qquad \mbox{for all } t\in [\tau,T-\tau]
  \eas
  and hence, by definition of $\I(T)$,
  \be{8.5}
	\I(T) \le C_7 + C_7 \I^\frac{3(p-\alpha)}{3p-1}(T)
	+ C_7\I^\frac{3(p-\alpha)q+2}{(3p-1)q}(T)
	\qquad \mbox{for all } T\in (2\tau,\tm-\tau)
  \ee
  with $C_7:=\frac{2C_5+C_6}{C_4}$.
  We can now rely on our assumptions that $\alpha>\frac{1}{3}$ and that (\ref{8.1}) holds, which namely ensure that
  \bas	
	\theta_1:=\frac{3(p-\alpha)}{3p-1}<\frac{3(p-\frac{1}{3})}{3p-1}=1
  \eas
  and
  \bas
	\theta_2:=\frac{3(p-\alpha)q+2}{(3p-1)q}
	= \frac{3(p-\alpha)+\frac{2}{q}}{3p-1}
	< \frac{3(p-\alpha)+(3\alpha-1)}{3p-1} =1,
  \eas
  respectively. Therefore, writing $\theta:=\max\{\theta_1, \theta_1\} \in (0,1)$ and noting that
  \bas
	\I(T) \le 2C_7 + 2C_7 \I^\theta(T)
	\qquad \mbox{for all } T\in (2\tau,\tm-\tau)
  \eas
  by (\ref{8.5}) and Young's inequality, we conclude by an elementary argument that
  \bas
	\I(T) \le C_8:=\max \Big\{ 1 \, , \, (4C_7)^\frac{1}{1-\theta}\Big\}
	\qquad \mbox{for all } T\in (2\tau,\tm-\tau).
  \eas
  In view of (\ref{8.4}), this in turn implies that
  \bas
	\int_t^{t+\tau} h(s) ds
	\le C_9:= C_5 + C_5 C_8^{\theta_1} + C_5 C_8^{\theta_2}
	\qquad \mbox{for all } t\in [\tau,\tm-\tau),
  \eas
  whereupon Lemma \ref{lem1} secondly guarantees that
  \bas
	y(t) \le C_9 + C_6
	\qquad \mbox{for all } t\in [\tau,\tm)
  \eas
  and thereby entails (\ref{8.2}), again because $n$ is bounded in $\Omega\times [0,\tau)$ by Lemma \ref{lem_loc}.
\qed
It remains to make sure 
that the above requirements on the auxiliary parameters $q$ and $r$ can indeed be fulfilled for arbitrarily large $p$
to end up with the following.
\begin{lem}\label{lem9}
  Suppose that $\alpha>\frac{1}{3}$. Then 
  given any $p>1$, one can find $C>0$ such that
  \be{9.1}
	\io n^p(\cdot,t) \le C
	\qquad \mbox{for all } t\in [0,\tm).
  \ee
\end{lem}
\proof
  As $\Omega$ is bounded, without loss of generality we may assume that $p$ additionally satisfies
  \be{9.11}
	p>\max\bigg\{ \frac{2}{3(3\alpha-1)} \, , \, \frac{1}{3\alpha} \, , \, \alpha+\frac{1}{3}\bigg\}.
  \ee
  We can then firstly pick $q>\frac{3}{2}$ such that
  \be{9.14}
	q<3p
  \ee
  and
  \be{9.2}
	q>\frac{2}{3\alpha-1}
  \ee
  as well as
  \be{9.3}
	q>\frac{1}{\alpha},
  \ee
  where the latter ensures that
  \bas
	\frac{(3p-1)q}{3}>\frac{3p-1}{3\alpha}.
  \eas
  Since furthermore our hypotheses that $p>1$ and $q>\frac{3}{2}$ warrant that
  \bas
	\frac{(3p-1)q}{3}>\frac{(3\cdot 1-1)\cdot\frac{3}{2}}{3}=1
  \eas
  and that clearly also $3p-1>1$, it is thereafter possible to choose $r>1$ in such a way that
  \be{9.4}
	\frac{3p-1}{3\alpha}<r<\min \Big\{ \frac{(3p-1)q}{3} \, , \, 3p-1\Big\}.
  \ee
  Now from (\ref{9.14}) and the third restriction in (\ref{9.11}) it follows that (\ref{3.1}) and (\ref{7.01}) hold,
  whereas (\ref{9.4}) guarantees validity of (\ref{5.1}), (\ref{6.2}) and (\ref{7.1}).
  As moreover (\ref{8.1}) is satisfied thanks to (\ref{9.2}), Lemma \ref{lem8} becomes applicable so as to assert the
  claimed boundedness property.
\qed
\mysection{Further regularity properties. Proof of Theorem \ref{theo14}}\label{sect_theo14}
Higher integrability properties can now be derived by applying arguments which are essentially standard in the
analysis of the heat and the Stokes equations.
Firstly, the uniform boundedness of $n$ with respect to the norm in $L^2(\Omega)$, together with our 
overall assumption that $f$ be bounded, entails the following.
\begin{lem}\label{lem10}
  Let $\alpha>\frac{1}{3}$. Then there exists $C>0$ such that
  \be{10.1}
	\|A^\beta u(\cdot,t)\|_{L^2(\Omega)} \le C
	\qquad \mbox{for all } t\in [0,\tm)
  \ee
  and
  \be{10.11}
	\|u(\cdot,t)\|_{L^\infty(\Omega)} \le C
	\qquad \mbox{for all } t\in [0,\tm).
  \ee
\end{lem}
\proof
  On the basis of a Duhamel formula associated with the Stokes subsystem of (\ref{0}), by means of well-known 
  smoothing properties of the Stokes semigroup (\cite{sohr}) we see that with some $\lambda_1>0$ and $C_1>0$ we have
  \bas
	\|A^\beta u(\cdot,t)\|_{L^2(\Omega)}
	&=& \bigg\| e^{-tA} A^\beta u_0 + \int_0^t A^\beta e^{-(t-s)A} 
		\proj\Big[ n(\cdot,s)\nabla \phi+f(\cdot,s)\Big] ds \bigg\|_{L^2(\Omega)} \\
	&\le& \|A^\beta u_0\|_{L^2(\Omega)} 
	+ C_1\int_0^t (t-s)^{-\beta} e^{-\lambda_1(t-s)} 
		\Big\| \proj\Big[n(\cdot,s)\nabla \phi+ f(\cdot,s)\Big]\Big\|_{L^2(\Omega)} ds
  \eas
  for all $t\in [0,\tm)$, because $\proj$ acts as an orthogonal projection on $L^2(\Omega;\R^3)$ (\cite{sohr}).
  Since Lemma \ref{lem9} together with the boundedness of $\nabla \phi$ and $f$ 
  entails the existence of $C_2>0$ such that $\|n(\cdot,s)\nabla \phi + f(\cdot,s)\|_{L^2(\Omega)} \le C_2$
  for all $s\in [0,\tm)$, and since $C_3:=\int_0^\infty \sigma^{-\beta} e^{-\lambda_1 \sigma} d\sigma$ is finite
  due to the fact that $\beta<1$, this implies that
  \bas
	\|A^\beta u(\cdot,t)\|_{L^2(\Omega)}
	\le \|A^\beta u_0\|_{L^2(\Omega)}
	+ C_1 C_2 C_3 \|\nabla\phi\|_{L^\infty(\Omega)}
	\qquad \mbox{for all } t\in [0,\tm)
  \eas
  and hence proves (\ref{10.1}), for $u_0\in D(A^\beta)$ by (\ref{init}).
  As our assumption $\beta>\frac{3}{4}$ warrants that $D(A^\beta) \hra L^\infty(\Omega;\R^3)$ 
  (\cite{giga1981_the_other}, \cite{henry}), this also entails (\ref{10.11}).
\qed
In conjunction again with Lemma \ref{lem9}, the latter entails a bound for $c$ in the flavor needed for an
application of Lemma \ref{lem_loc} for the derivation of Theorem \ref{theo14}.
\begin{lem}\label{lem11}
  If $\alpha>\frac{1}{3}$, then for all $p>1$ there exists $C>0$ such that
  \be{11.1}
	\|c(\cdot,t)\|_{W^{1,p}(\Omega)} \le C
	\qquad \mbox{for all } t\in [0,\tm).
  \ee
\end{lem}
\proof 
  We let $B$ denote the realization of $-\Delta+\frac{1}{2}$ under homogeneous Neumann boundary conditions 
  in $L^p(\Omega)$ and then obtain that $B$ is sectorial with its spectrum contained in $[\frac{1}{2},\infty)$,
  and that for arbitrary $\gamma\in (\frac{1}{2},1)$ the corresponding fractional power $B^\gamma$ has its domain
  satisfy $D(B^\gamma) \hra W^{1,p}(\Omega)$ (\cite{henry}), so that
  \be{11.2}
	\|\varphi\|_{W^{1,p}(\Omega)} \le C_1(\gamma) \|B^\gamma \varphi\|_{L^p(\Omega)}
	\qquad \mbox{for all } \varphi\in D(B^\gamma)
  \ee
  with some $C_1(\gamma)>0$.
  Hencoforth fixing any $\gamma\in (\frac{1}{2},1)$ and $\gamma'\in (\frac{1}{2},\gamma)$, by a well-known 
  interpolation property (\cite{friedman}) we can find $C_2>0$ such that
  \be{11.3}
	\|B^{\gamma'}\varphi\|_{L^p(\Omega)}
	\le C_2 \|B^\gamma\varphi\|_{L^p(\Omega)}^a \|\varphi\|_{L^p(\Omega)}^{1-a}
	\qquad \mbox{for all } \varphi\in D(B^\gamma)
  \ee
  with $a:=\frac{\gamma'}{\gamma}\in (0,1)$, where according to the Gagliardo-Nirenberg inequality there exists $C_3>0$
  fulfilling
  \bas
	\|\varphi\|_{L^p(\Omega)}^{1-a} 
	\le C_3\|\varphi\|_{W^{1,p}(\Omega)}^{(1-a)b} \|\varphi\|_{L^1(\Omega)}^{(1-a)(1-b)} 
	\qquad \mbox{for all } \varphi\in W^{1,p}(\Omega)
  \eas
  with $b:=\frac{3(p-1)}{4p-3}\in (0,1)$. Together with (\ref{11.2}) and (\ref{11.3}), this shows that if we let
  $d:=a+(1-a)b\in (0,1)$ and $C_4:=C_1^{(1-a)b}(\gamma) C_2 C_3$, then
  \be{11.4}
	\|B^{\gamma'}\varphi\|_{L^p(\Omega)}
	\le C_4\|B^\gamma\varphi\|_{L^p(\Omega)}^d \|\varphi\|_{L^1(\Omega)}^{1-d}
	\qquad \mbox{for all } \varphi\in D(B^\gamma).
  \ee
  Now since $c_t=-(B+\frac{1}{2})c+n-u\cdot\nabla c$ in $\Omega\times (0,\tm)$ by (\ref{0}), an associated
  variation-of-constants representation together with known regularization features of the corresponding
  analytic semigroup $(e^{-tB})_{t\ge 0}$ shows that there exists $C_5>0$ such that
  \bea{11.5}
	\hspace*{-10mm}
	\|B^\gamma c(\cdot,t)\|_{L^p(\Omega)}
	&=& \bigg\| B^\gamma e^{-(t-\tau)(B+\frac{1}{2})} c(\cdot,\tau)
	+ \int_\tau^t B^\gamma e^{-(t-s)(B+\frac{1}{2})} n(\cdot,s) ds \nn\\
	& & - \int_\tau^t B^\gamma e^{-(t-s)(B+\frac{1}{2})} u(\cdot,s)\cdot\nabla c(\cdot,s) ds
	\bigg\|_{L^p(\Omega)} \nn\\
	&\le& C_5\|B^\gamma c(\cdot,\tau)\|_{L^p(\Omega)}
	+ C_5\int_\tau^t (t-s)^{-\gamma} e^{-\frac{1}{2}(t-s)} \|n(\cdot,s)\|_{L^p(\Omega)} ds \nn\\
	& & + C_5\int_\tau^t (t-s)^{-\gamma} e^{-\frac{1}{2}(t-s)} 
		\Big\| u(\cdot,s)\cdot\nabla c(\cdot,s)\Big\|_{L^p(\Omega)} ds 
	\qquad \mbox{for all } t\in [\tau,\tm).
  \eea
  Here by Lemma \ref{lem9} we can find $C_6>0$ such that
  \be{11.6}
	\|n(\cdot,s)\|_{L^p(\Omega)} \le C_6
	\qquad \mbox{for all } s\in [\tau,\tm),
  \ee
  while Lemma \ref{lem10} together with (\ref{11.2}), (\ref{11.4}) and (\ref{mass_c}) shows that with some
  $C_7>0$ we have
  \bas
	\Big\|u(\cdot,s)\cdot\nabla c(\cdot,s)\Big\|_{L^p(\Omega)}
	&\le& \|u(\cdot,s)\|_{L^\infty(\Omega)} \|\nabla c(\cdot,s)\|_{L^p(\Omega)} \\
	&\le& C_7\|\nabla c(\cdot,s)\|_{L^p(\Omega)} \\
	&\le& C_1(\gamma') C_7\|B^{\gamma'} c(\cdot,s)\|_{L^p(\Omega)} \\
	&\le& C_1(\gamma') C_4 C_7 \|B^\gamma c(\cdot,s)\|_{L^p(\Omega)}^d \|c(\cdot,s)\|_{L^1(\Omega)}^{1-d} \\
	&\le& C_8 \|B^\gamma c(\cdot,s)\|_{L^p(\Omega)}^d
	\qquad \mbox{for all } s\in [\tau,\tm)
  \eas
  with $C_8:=C_1(\gamma') C_2 C_5 \cdot \max \Big\{ \io c_0, \io n_0\Big\}$.
  Combining this with (\ref{11.6}) and (\ref{11.5}) and abbreviating 
  $C_9:=\int_0^\infty \sigma^{-\gamma} e^{-\frac{\sigma}{2}} d\sigma<\infty$ as well as 
  $M(T):=\sup_{t\in [\tau,T]} \|B^\gamma c(\cdot,t)\|_{L^p(\Omega)}$ for $T\in (\tau,\tm)$, we obtain
  \bas
	\|B^\gamma c(\cdot,t)\|_{L^p(\Omega)}
	\le C_5 \|B^\gamma c(\cdot,\tau)\|_{L^p(\Omega)}
	+ C_5 C_6 C_9 + C_5 C_8 C_9 M^d(T)
	\qquad \mbox{for all } t\in [\tau,T]
  \eas
  and hence
  \bas
	M(T) \le C_{10} + C_{10} M^d(T)
	\qquad \mbox{for all } T\in [\tau,\tm)
  \eas
  with $C_{10}:=\max \big\{ C_5\|B^\gamma c(\cdot,\tau)\|_{L^p(\Omega)} + C_5 C_6 C_9, C_5 C_8 C_9\big\}$.
  As $d<1$, this entails that $M(T)\le \max\{1,(2C_{10})^\frac{1}{1-d}\}$ for all $T\in [\tau,\tm)$, which in light of
  (\ref{11.2}) establishes (\ref{11.1}) due to the inclusion $c\in L^\infty((0,\tau);W^{1,p}(\Omega))$ asserted by
  Lemma \ref{lem_loc}.
\qed
Finally, pointwise boundednes of $n$ results from a standard argument contained in the literature.
\begin{lem}\label{lem12}
  Let $\alpha>\frac{1}{3}$. Then with some $C>0$ we have
  \be{12.1}
	\|n(\cdot,t)\|_{L^\infty(\Omega)} \le C
	\qquad \mbox{for all } t\in [0,\tm).
  \ee
\end{lem}
\proof
  We write the first equation in (\ref{0}) in the form $n_t=\Delta n + \nabla \cdot h(x,t)$ 
  with $h:=-nS(n)\nabla c-nu$ and then obtain from (\ref{0}) that $h\cdot\nu=0$
  on $\pO\times (0,\tm)$, whereas (\ref{S}), Lemma \ref{lem9}, Lemma \ref{lem11} and Lemma \ref{lem10} entail that
  $h\in L^\infty((0,\tm);L^p(\Omega;\R^3))$ for each $p>1$.
  Since moreover $n\in L^\infty((0,\tm);L^p(\Omega))$ for any such $p$, (\ref{12.1}) can e.g.~be derived by
  a Moser-type iterative argument; for a statement precisely covering the present situation we refer to
  \cite[Lemma A.1]{taowin_subcrit}.
\qed
Thanks to the extensibility criterion (\ref{ext}), 
the derivation of our main results thus consists in a mere collection of the latter three lemmata.\abs
\proofc of Theorem \ref{theo14}. \quad
  In view of Lemma \ref{lem_loc}, the boundedness properties obtained Lemma \ref{lem12}, Lemma \ref{lem11}
  and Lemma \ref{lem10} assert both global extensibility and the claimed regularity features
  of the local-in-time solution from Lemma \ref{lem_loc},
  as well as the temporally uniform estimate in (\ref{14.2}).
\qed
\vspace*{5mm}
{\bf Acknowledgement.} \quad
  The author acknowledges support of the {\em Deutsche Forschungsgemeinschaft} in the context of the project
  {\em Analysis of chemotactic cross-diffusion in complex frameworks}, and he is grateful to Yulan Wang for
  numerous helpful remarks on this manuscript.
\end{document}